# Independent transversal blow-up of graphs

Tianjiao Dai[1]   Weichan Liu[1*]   Xin Zhang[2]

1. School of Mathematics, Shandong University, Jinan, 250100, China

tianjiao.dai@sdu.edu.cn, wcliu@sdu.edu.cn

2. School of Mathematics and Statistics, Xidian University, Xi'an, 710071, China

xzhang@xidian.edu.cn

February 27, 2025

**Abstract**

In an $r$-partite graph, an independent transversal of size $s$ (ITS) consists of $s$ vertices from each part forming an independent set. Motivated by a question from Bollobás, Erdős, and Szemerédi (1975), Di Braccio and Illingworth (2024) inquired about the minimum degree needed to ensure an $n \times \cdots \times n$ $r$-partite graph contains $K_r(s)$, a complete $r$-partite graph with $s$ vertices in each part. We reformulate this as finding the smallest $n$ such that any $n \times \cdots \times n$ $r$-partite graph with maximum degree $\Delta$ has an ITS. For any $\varepsilon > 0$, we prove the existence of a $\gamma > 0$ ensuring that if $G$ is a multipartite graph partitioned as $(V_1, V_2, \ldots, V_r)$, where the average degree of each part $V_i$ is at most $D$, the maximum degree of any vertex to any part $V_i$ is at most $\gamma D$, and the size of each part $V_i$ is at least $(s + \varepsilon)D$, then $G$ possesses an ITS. The constraint $(s + \varepsilon)D$ on the part size is tight. This extends results of Loh and Sudakov (2007), Glock and Sudakov (2022), and Kang and Kelly (2022). We also show that any $n \times \cdots \times n$ $r$-partite graph with minimum degree at least $\left(r - 1 - \frac{1}{2s^2}\right) n$ contains $K_r(s)$ and provide a relative Turán-type result. Additionally, this paper explores counting ITSs in multipartite graphs.

# 1 Introduction

## 1.1 Background and Research Motivation

Turán-type problems occupy a significant and prominent position within the broader realm of Extremal Graph Theory, serving as pivotal areas of research that explore the boundaries of graph properties and their associated extremal values. Given a graph $F$, we define $\mathrm{ex}(n, F)$ as the maximum number of edges that a graph on $n$ vertices can have without containing $F$ as a subgraph. We denote by $K_{s_1,\dots,s_r}$ the complete $r$-partite graph where the parts have sizes $s_1, \dots, s_r$, and we use the shorthand $K_r(s)$ to represent $K_{s,\dots,s}$ with $r$ parts of equal size $s$.

*Corresponding author.

A renowned theorem by Erdős and Stone [8] establishes the asymptotic behavior of $\mathrm{ex}(n, K_{r+1}(s))$:

$$\mathrm{ex}(n, K_{r+1}(s)) = \mathrm{ex}(n, K_{r+1}) + o(n^2) = \left(1 - \frac{1}{r}\right)\frac{n^2}{2} + o(n^2).$$

Subsequently, Erdős [7] and Simonovits [23] independently refined the error term to $O(n^{2-\frac{1}{s}})$.

Since 1951, when Zarankiewicz [27] suggested exploring the maximum number of edges in a $K_{s,t}$-free bipartite graph, extremal problems pertaining to multipartite graphs have garnered significant attention and study. In 1975, Bollobás, Erdős, and Szemerédi [3] introduced the following problem which is a tripartite version of the Erdős-Stone theorem / Zarankiewicz problem. It's worth noting that they opted to frame this problem in terms of minimum degree instead of the number of edges, considering that it exists in a sparser context compared to the edge version.

**Problem 1.** [3, Bollobás, Erdős, and Szemerédi] *For each $s$, what is the smallest $\tau = \tau(n)$ such that any $n \times n \times n$ tripartite graph with minimum degree at least $n + \tau$ contains $K_3(s)$?*

In addressing this problem, Bollobás, Erdős, and Szemerédi [3] showed that $\tau = O(n^{1-\frac{1}{3s^2}})$ and conjectured $\tau = O(n^{1/2})$ when $s = 2$. This conjecture was very recently confirmed by Di Braccio and Illingworth [6], and they also improved the answer to Problem 1 from $O(n^{1-\frac{1}{3s^2}})$ to $O(n^{1-\frac{1}{s}})$. Note that recently, Chen et al. [5] also demonstrated a weaker result, specifically that $\tau = O(n^{1-\frac{1}{s^2+s}})$.

Let $G_r(n)$ be an $r$-partite graph with $n$ vertices in each block (here a block refers to one of the $r$ disjoint partitions that the vertex set of the graph is divided into) and let $f(n) = \max\{\delta(G) : K_r \nsubseteq G_r(n)\}$. Set $c_r = \lim_{n\to\infty} \frac{f_r(n)}{n}$. Bollobás, Erdős, and Szemerédi [3] proved that

$$\lim_{r\to\infty}(c_r - (r-2)) \geq \frac{1}{2}$$

and conjectured that equality holds. This conjecture was later confirmed by Haxell [13].

Exploring an $r$-partite variant of Problem 1 now appears to be a natural extension. Di Braccio and Illingworth, in their recent work [6], posed the following question:

**Problem 2.** [6, Di Braccio and Illingworth] *For each $s$, what is the smallest $\tau = \tau(n)$ such that any $n \times \cdots \times n$ $r$-partite graph with minimum degree at least $c_r n + \tau$ contains $K_r(s)$?*

Problem 2 with $s = 1$ ignited a sequence of research [13, 15], with Szabó and Tardos [24] providing the solution for even $r$, and Haxell and Szabó [12] offering the solution for odd $r$. Moreover, the multi-partite version of the Erdős-Stone theorem was recently investigated by Han and Zhao [11]. Indeed, they studied $\mathrm{ex}_r(n, K_r(s))$, which defines the largest number of edges in $K_r(s)$-free graphs from the family of $r$-partite graphs with $n$ vertices in each part.

Given an $r$-partite graph $G$ with a partition $\mathcal{P} = (V_1, V_2, \ldots, V_r)$ such that each $V_i$ is an independent set, we construct another $r$-partite graph $\tilde{G}$ with the same partition $\mathcal{P}$ (we call $\tilde{G}$ an $r$-partite complement of $G$). The additional condition for $\tilde{G}$ is that for any two vertices $u$ and $v$ from different parts $V_i$ and $V_j$ (where $i \neq j$), the pair $uv$ forms an edge in $\tilde{G}$ if and only if $uv$ does not form an edge in $G$. Now, if $G$ is an $n \times \cdots \times n$ $r$-partite graph with minimum degree at least $n + \tau$, then $\tilde{G}$ is an $n \times \cdots \times n$ $r$-partite graph with maximum degree at most $(r-1)n - (c_r n + \tau) = (r-1-c_r)n - \tau$. Alternatively, a $K_r(s)$ structure in $G$ maps to an independent set in $\tilde{G}$, composed of precisely $s$ vertices drawn from each subset $V_i$. In technical terms, we define this type of independent set as an *independent transversal blow-up* of size $s$ in $\tilde{G}$ with respect to the partition $\mathcal{P}$. In the special case where $s$ equals 1, it is simply designated as an *independent transversal*. To clarify, in the context of an $r$-partite graph, a transversal consists of a single vertex from each block, whereas an independent transversal is both an independent set and a transversal. Leveraging this transformation, we restate Problem 2 as follows:

**Problem 3.** *For each $s$, what is the smallest $\tau = \tau(n)$ such that any $n \times \cdots \times n$ $r$-partite graph with maximum degree at most $(r-1-c_r)n - \tau$ contains an independent transversal blow-up of size $s$?*

By substituting $\Delta := (r-1-c_r)n - \tau$ into Problem 3, we notice that minimizing $\tau$ is equivalent to determining the smallest block size, given the relationship $n = \frac{\Delta+\tau}{r-1-c_r}$. Therefore, we reformulate Problem 3 as follows:

**Problem 4.** *For each $s$, what is the smallest $n = n(\Delta)$ such that any $n \times \cdots \times n$ $r$-partite graph with maximum degree at most $\Delta$ contains an independent transversal blow-up of size $s$?*

Problem 4 with $s = 1$ sparked a series of research endeavors. Alon [2] used the Lovász Local Lemma to show that $n \geq 2e\Delta$. Haxell [13] later improved this result to $n \geq 2\Delta$. Szabó and Tardos [24] constructed graphs with parts of size $2\Delta - 1$ and maximum degree $\Delta$, with no independent transversals, therefore, the bound $2\Delta$ is tight. This subsequently provides a solution of $\tau \leq (r-2.5)n$ for both Problem 2 and Problem 3 when $s = 1$.

In the construction of Szabó and Tardos [24], the graphs are disjoint unions of $2\Delta - 1$ complete bipartite subgraphs $K_{\Delta,\Delta}$, and the partition into $V_1, V_2, \ldots, V_r$ is done in such a way that the parts $V_1, V_2, \ldots, V_r$ separate the sides of each $K_{\Delta,\Delta}$ (thus $r$ is $2\Delta$ in this construction). Ideally, by setting $r$ to a constant that is independent of $\Delta$, we could potentially enhance the bound for $\tau$ from $O(n)$ to $O(n^{1-\varepsilon})$ for some $\varepsilon > 0$. However, we will not explore this direction further in this paper.

## 1.2 Locally sparseness

In 2007, Loh and Sudakov [18] proposed another direction to reduce the bound of $2\Delta$ in Haxell's theorem [13]. They introduced a constraint on the local degree. Formally, given an $r$-partite graph $G$ with a partition $\mathcal{P} = (V_1, V_2, \ldots, V_r)$, we define the *local degree* of $G$ with respect to $\mathcal{P}$ as

$$\Gamma_{\mathcal{P}}(G) := \max_{v \in V} \max_{i \in [r]} \{d_G(v, V_i)\}.$$

In the construction by Szabó and Tardos [24], numerous pairs of disjoint sets $(V_i, V_j)$ exist, with complete bipartite subgraphs of linear size connecting them. The number of edges between such pairs $(V_i, V_j)$ scales quadratically with $\Delta$, resulting in a local degree of $O(\Delta)$ in this construction.

This local degree constraint naturally arises in various situations, with vertex list coloring being one such example. Given a graph $G$ with a list assignment $L(v)$ of colors to each vertex $v \in V(G)$, our task is to color the vertices so that each vertex receives a color from its list and no two adjacent vertices share the same color. Additionally, assume that for each vertex $v$, any color $c$ appears in the lists of at most $\Delta$ of its neighboring vertices. Our primary concern is identifying the smallest possible size of the lists that guarantees a proper coloring. Reed and Sudakov [22] demonstrated that lists of size $(1 + o(1))\Delta$ are actually sufficient.

This issue closely mirrors Problem 4 when $s$ is set to 1. In fact, we expand each vertex $v \in V(G)$ into $|L(v)|$ vertices, denoted as $(v, c)$ for each $c \in L(v)$. We connect two vertices $(v, c)$ and $(v', c')$ with an edge if $vv' \in E(G)$ and $c = c'$. This transformation leads to an $r$-partite graph $H$ (where $r$ equals the number of vertices in $G$) with a maximum degree of $\Delta$ and a local degree of 1. Consequently, an independent transversal in $H$ corresponds directly to a proper list coloring of $G$. Therefore, the result presented by Reed and Sudakov [22] suggests that every $r$-partite graph with a maximum degree of $\Delta$ and a local degree of 1 possesses an independent transversal. Subsequently, Loh and Sudakov [18] reached the same conclusion by modifying the condition to a local degree of $o(\Delta)$.

In 2022, Glock and Sudakov [10], and independently, Kang and Kelly [16], substituted the maximum degree condition from Loh and Sudakov's result [18] with an average degree condition, significantly enhancing the applicability of the outcome. Given an $r$-partite graph $G$ with a partition $\mathcal{P} = (V_1, V_2, \ldots, V_r)$ such that each $V_i$ is an independent set for each $i \in [r]$, we define the *average degree* of $V_i$ in $G$ as:

$$\overline{d}_G(V_i) := \frac{1}{|V_i|} \sum_{v \in V_i} d_G(v).$$

Additionally, we define the *maximum block average degree* as:

$$\overline{D}_{\mathcal{P}}(G) := \max_{i \in [r]} \overline{d}_G(V_i)$$

This represents the highest average degree among all the independent sets $V_i$ in the partition $\mathcal{P}$. Both research groups demonstrated that every $r$-partite graph with a maximum block average degree of $D$ and a local degree of $o(D)$ admits an independent transversal.

## 1.3 Main Results

In this paper, we mainly proved Theorem 4.8 which extends the findings of Loh and Sudakov [18], Glock and Sudakov [10], and Kang and Kelly [16]. Using a semi-random algorithm approach, we integrate ideas from three papers with novel innovations. Specifically, we introduce two auxiliary graphs in Section 3 linking the construction of independent transversals and the process of independent transversal blow-ups, strategically designed to facilitate probabilistic methods in proving the main theorem while reducing dependence risks.

It is worth noting that simplifying the block size condition from $(s+\varepsilon)D$ to $(s^2+\varepsilon)D$ would make the proof quite easier, see Theorem 4.7. Subsequently, we employ Theorem 4.17 to demonstrate that the block size condition $(s+\varepsilon)D$ is optimal. A partition $\mathcal{P} = (V_1, V_2, \ldots, V_r)$ is *$s$-thick* if each $V_i$ has size at least $s$.

**Theorem 4.8.** *For every $\varepsilon > 0$, there exists $\gamma > 0$ such that the following hold. If $G$ is a multipartite graph with a partition $\mathcal{P} = (V_1, V_2, \ldots, V_r)$, where $\overline{D}_{\mathcal{P}}(G) \leq D$, $\Gamma_{\mathcal{P}}(G) \leq \gamma D$, and $\mathcal{P}$ is $(s+\varepsilon)D$-thick, then $G$ has an independent transversal blow-up with size $s$ with respect to $\mathcal{P}$.*

**Theorem 4.17.** *There is a multipartite graph $G$ and a $sD$-thick partition $\mathcal{P}$ such that the maximum block average degree is $D$, the local degree is 1, and there is no independent transversal blow-up with size $s$ with respect to $\mathcal{P}$.*

Disregarding the local degree condition, we obtain two theorems: one for the maximum block average degree setting and one for the maximum degree setting.

**Theorem 4.3.** *Given a multipartite graph $G$ with a partition $\mathcal{P} = (V_1, V_2, \ldots, V_r)$, if $\overline{D}_{\mathcal{P}}(G) \leq D$ and $\mathcal{P}$ is $4s^2D$-thick, then $G$ has an independent transversal blow-up with size $s$ with respect to $\mathcal{P}$.*

**Theorem 4.5.** *Given a multipartite graph $G$ with a partition $\mathcal{P} = (V_1, V_2, \ldots, V_r)$, if $\Delta(G) \leq \Delta$ and $\mathcal{P}$ is $2s^2\Delta$-thick, then $G$ has an independent transversal blow-up with size $s$ with respect to $\mathcal{P}$.*

Set $n = 2s^2\Delta$ and let $G$ be an $n\times\cdots\times n$ $r$-partite graph with minimum degree at least $(r-1-\frac{1}{2s^2})n$. Now, the $r$-partite complement $\tilde{G}$ of $G$ has maximum degree at most $(r-1)n-(r-1-\frac{1}{2s^2})n = \Delta$ and thus $G$ has an independent transversal blow-up with size $s$ with respect to the partition by Theorem 4.5. This in turn indicates that $G$ contains $K_r(s)$, thereby providing a solution to Problems 2 and 3 with the conclusion that

$$\tau \leq \left(r - 2 - \frac{1}{2s^2}\right)n.$$

Next, set $n = 4s^2D$ and let $G$ be an $n\times\cdots\times n$ $r$-partite graph such that $\sum_{v\in V_i} d_G(v) \geq (r-1-\frac{1}{4s^2})n^2$ for each part $V_i$. In the $r$-partite complement $\tilde{G}$ of $G$, we have $\sum_{v\in V_i} d_{\tilde{G}}(v) = \sum_{v\in V_i}((r-1)n-d_G(v)) \leq (r-1)n^2 - (r-1-\frac{1}{4s^2})n^2 = nD$. By applying Theorem 4.3, we have the following Turán-type answer to Problem 2..

**Corollary 1.1.** *Any $n\times\cdots\times n$ $r$-partite graph, where each part is incident with at least $\left(r-1-\frac{1}{4s^2}\right)n^2$ edges, necessarily contains $K_r(s)$.*

A *factor of transversal* is a set of pairwise-disjoint transversals covering all vertices. The problem of finding sufficient conditions for the existence of factors of independent transversals in multipartite graphs was investigated by several researchers [1, 4, 18]. This question is also related to the notion of strong chromatic number [9, 14]. In this paper, we prove the following two theorems on the factor of independent transversal blow-ups.

**Theorem 5.4.** *Given a multipartite graph $G$ with a partition $\mathcal{P} = (V_1, V_2, \ldots, V_r)$ such that $|V_i| = 3s^2\Delta$ for each $i$, if $\Delta(G) \leq \Delta$, then $G$ has a factor of independent transversal blow-ups with size $s$ with respect to $\mathcal{P}$.*

**Theorem 5.6.** *For every $\varepsilon > 0$, the following holds. Given a multipartite graph $G$ with a partition $\mathcal{P} = (V_1, V_2, \ldots, V_r)$ such that $|V_i| = (2 + \varepsilon)s^2\Delta$ for each $i$, if $\Delta(G) \leq \Delta$ and $\Gamma_{\mathcal{P}}(G) = o(\Delta)$, then $G$ has a factor of independent transversal blow-ups with size $s$ with respect to $\mathcal{P}$.*

The final contribution of this paper is to count the number of independent transversal blow-ups, without requiring that any two of them are vertex-disjoint. The theorem below states that if $G$ is a multipartite graph with a maximum block average degree of at most $D$ and a block size of at least $4s^2D$, then the number of such independent transversal blow-ups of size $s$ contained in $G$ is approximately $\Omega(D^r)$.

**Theorem 5.8.** *Given a multipartite graph $G$ with a partition $\mathcal{P} = (V_1, V_2, \ldots, V_r)$, if for each $i$, $|V_i| \geq t$ and $V_i$ is incident with at most $\frac{t}{4s^2}|V_i|$ edges, then $G$ has at least*

$$\frac{1}{2^r}\binom{t}{s}^r$$

*independent transversal blow-ups with size $s$ with respect to $\mathcal{P}$. Consequently, if $\Gamma_{\mathcal{P}}(G) \leq D$ and $\mathcal{P}$ is $4s^2D$-thick, then then $G$ possesses at least as many independent transversal blow-ups of size $s$ with respect to $\mathcal{P}$ as there are when substituting $4s^2D$ for $t$.*

## 1.4 Organization of the paper

The paper is organized as follows:

- **Section 2**: We present the probabilistic tools that are crucial for our proofs.
- **Section 3**: We introduce two types of auxiliary graphs and explore their properties, which will be essential for our subsequent analysis.
- **Section 4**: Our main focus is on demonstrating the presence of an independent transversal blow-up in multipartite graphs with specific properties. Specifically, our main theorem (Theorem 4.8) and Theorem 4.17 will be rigorously proven in Subsection 4.2. Additionally, other theorems, including Theorems 4.3 and 4.5, will be presented in Subsection 4.1.

- **Section 5**: We demonstrate the existence of a factor of independent transversal blow-ups in multipartite graphs with specific properties using Theorems 5.4 and 5.6, and also count how many of these independent transversal blow-ups exist in multipartite graphs with different properties, with the help of Theorem 5.8.

- **Section 6**: We provide comments on our results and propose several intriguing problems for future investigation.

# 2 Probabilistic tools

**Lemma 2.1.** [19, Chernoff's Bound] *Let $X$ be the sum of independent Bernoulli random variables, where each may have a distinct expected value. For any $0 \le s \le \mathbb{E}[X]$, the following bound holds:*

$$\mathbb{P}[|X - \mathbb{E}[X]| \ge s] \le 2\exp\left(-\frac{s^2}{3\mathbb{E}[X]}\right).$$

*Additionally, for $s \ge 7\mathbb{E}[X]$, the following bound is given:*

$$\mathbb{P}[X \ge s] \le e^{-s}.$$

**Lemma 2.2.** [21, Talagrand's Inequality] *Consider a non-negative random variable $X$ that depends on $n$ independent trials $T_1, \dots, T_n$. Suppose $X$ satisfies the following properties:*

1. ***C-Lipschitz**: Changing the outcome of any single trial can change $X$ by at most $C$.*

2. ***r-Certifiable**: If $X \ge s$ for some $s > 0$, then there exists a set of at least $rs$ trials whose outcomes confirm that $X \ge s$.*

*Then, for any $t \ge 0$, the following inequality holds:*

$$\mathbb{P}[|X - \mathbb{E}[X]| > t + 20C\sqrt{r\mathbb{E}[X]} + 64C^2 r] \le 4\exp\left(-\frac{t^2}{8C^2 r(\mathbb{E}[X] + t)}\right).$$

**Remark:** This inequality does not appear in its original form in Talagrand's paper [25]. However, for the sake of convenience in combinatorial applications, it has undergone numerous variations, including those by Reed and Sudakov [22], Loh and Sudakov [18], Glock and Sudakov [10], Molloy and Reed [21], Kelly and Postle [17], and Chakraborti and Tran [4].

**Lemma 2.3.** [20, Lovász Local Lemma] *Let $\mathcal{B}$ be a collection of bad events. Suppose that:*

1. *Each event $B$ in $\mathcal{B}$ occurs with probability at most $p$.*

2. *Each event $B$ in $\mathcal{B}$ is independent of all but at most $d$ other events in $\mathcal{B}$.*

*If $ep(d + 1) < 1$, then there is a positive probability that none of the bad events in $\mathcal{B}$ will occur.*

## 3 Auxiliary graphs

**Definition 1.** *Let $G$ be a multipartite graph with a partition $\mathcal{P} = (V_1, V_2, \ldots, V_r)$. Suppose each subset $V_i$ has size $\ell > s$. We define an auxiliary graph $\dot{G}$ based on $G$ and $\mathcal{P}$ as follows:*

1. *For each integer $i$ and each $s$-set $Z \subseteq V_i$ in $G$, we create a vertex $v_Z$ in $\dot{G}$.*

2. *Two vertices $v_{Z_1}$ and $v_{Z_2}$ in $\dot{G}$ are adjacent if and only if there exists at least one edge $e \in E(G)$ such that $e$ connects a vertex in $Z_1$ to a vertex in $Z_2$.*

*It is evident that $\dot{G}$ is a multipartite graph with a partition $\dot{\mathcal{P}}$ comprising $\bigcup_i \dot{V}_i$, where each $\dot{V}_i$ encompasses all vertices $v_Z$ in $\dot{G}$ such that $Z$ is an $s$-set drawn from $V_i$ in $G$.*

**Observation 1.** *The following holds.*

1. $|\dot{V}_i| = \binom{\ell}{s}$
2. $\Delta(\dot{G}) \leq s\binom{\ell-1}{s-1}\Delta(G)$
3. $\Gamma_{\dot{\mathcal{P}}}(\dot{G}) \leq s\binom{\ell-1}{s-1}\Gamma_{\mathcal{P}}(G)$
4. $\overline{D}_{\dot{\mathcal{P}}}(\dot{G}) \leq s\binom{\ell-1}{s-1}\overline{D}_{\mathcal{P}}(G)$

**Proof.** The first item is straightforward so we prove the remaining ones. Fix a vertex $v_{Z_1} \in \dot{V}_i$, where $Z_1 = \{v_1, v_2, \ldots, v_s\}$. For each neighbor $u$ of $v_1$, $u$ is contained in $\binom{\ell-1}{s-1}$ $s$-sets in its part. This implies that $v_{Z_1}$ has $\binom{\ell-1}{s-1}$ neighbors in $\dot{G}$ due to the edge $v_1 u$ in $G$.

Since $G$ has at most $s\Delta(G)$ edges in the form $v_q u$ with $v_q \in Z_1$, $v_{Z_1}$ has at most $s\binom{\ell-1}{s-1}\Delta(G)$ neighbors in $\dot{G}$. Since each $v_q$ has at most $\Gamma_{\mathcal{P}}(G)$ neighbors in any block $V_j$, a similar argument implies that $v_{Z_1}$ has at most $s\binom{\ell-1}{s-1}\Gamma_{\mathcal{P}}(G)$ neighbors in $\dot{V}_j$ for any $j \neq i$. For the average degree,

$$\begin{aligned}\sum_{v_Z \in \dot{V}_i} d_{\dot{G}}(v_Z) &\leq \sum_{Z \subset \binom{[\ell]}{s}} \sum_{v \in Z} \binom{\ell-1}{s-1} d_G(v) = \binom{\ell-1}{s-1} \sum_{Z \subset \binom{V_i}{s}} \sum_{v \in Z} d_G(v) \\ &= \binom{\ell-1}{s-1}\binom{\ell-1}{s-1} \sum_{v \in V_i} d_G(v) \leq \binom{\ell-1}{s-1}^2 \ell \overline{D}_{\mathcal{P}}(G),\end{aligned}$$

implying

$$\overline{d}_{\dot{G}}(\dot{V}_i) = \frac{\sum_{v_Z \in \dot{V}_i} d_{\dot{G}}(v_Z)}{|\dot{V}_i|} \leq \frac{\binom{\ell-1}{s-1}^2 \ell \overline{D}_{\mathcal{P}}(G)}{\binom{\ell}{s}} = s\binom{\ell-1}{s-1}\overline{D}_{\mathcal{P}}(G).$$

This completes the proof of items 2–4. □

**Definition 2.** *Let $G$ be a multipartite graph with a partition $\mathcal{P} = (V_1, V_2, \ldots, V_r)$. Suppose each subset $V_i$ can be further subdivided into $t$ subsets of size $s$, denoted as $V_i^1, V_i^2, \ldots, V_i^t$. We define an auxiliary graph $\ddot{G}$ based on $G$ and this refined partition $\mathcal{P}$ as follows:*

1. *For each integer $i$ and each $s$-set $Z \in \{V_i^1, V_i^2, \ldots, V_i^t\}$ in $G$, we create a vertex $v_Z$ in $\ddot{G}$.*

2. *Two vertices $v_{Z_1}$ and $v_{Z_2}$ in $\ddot{G}$ are adjacent if and only if there exists at least one edge $e \in E(G)$ such that $e$ connects a vertex in $Z_1$ to a vertex in $Z_2$.*

*It is clear that $\ddot{G}$ is a multipartite graph with a partition $\ddot{\mathcal{P}}$ consisting of the union $\bigcup_i \ddot{V}_i$, where each $\ddot{V}_i$ contains all vertices $v_Z$ in $\ddot{G}$ for which $Z$ belongs to the set $\{V_i^1, V_i^2, \dots, V_i^t\}$.*

**Observation 2.** *The following holds.*

1. $|\ddot{V}_i| = t$
2. $\Delta(\ddot{G}) \le s\Delta(G)$
3. $\Gamma_{\ddot{\mathcal{P}}}(\ddot{G}) \le s\Gamma_{\mathcal{P}}(G)$
4. $\overline{D}_{\ddot{\mathcal{P}}}(\ddot{G}) \le s\overline{D}_{\mathcal{P}}(G)$

**Proof.** The proofs for the first three items are omitted here as they are straightforward. For the last one, since

$$\sum_{v_Z \in \ddot{V}_i} d_{\ddot{G}}(v_Z) \le \sum_{j=1}^{t} \sum_{v \in V_i^j} d_G(v) = \sum_{v \in V_i} d_G(v) \le ts\overline{D}_{\mathcal{P}}(G),$$

we have

$$\overline{d}_{\ddot{G}}(\ddot{V}_i) = \frac{\sum_{v_Z \in \ddot{V}_i} d_{\ddot{G}}(v_Z)}{|\ddot{V}_i|} \le \frac{ts\overline{D}_{\mathcal{P}}(G)}{t} = s\overline{D}_{\mathcal{P}}(G),$$

as desired. □

# 4 One independent transversal blow-up

## 4.1 Warm-up

**Theorem 4.1.** *Given a multipartite graph $G$ with a partition $\mathcal{P} = (V_1, V_2, \dots, V_r)$, if $\overline{d}_G(V_i) \le D$ for each $i \in [r]$ and $\mathcal{P}$ is $2es^2D$-thick, then $G$ has an independent transversal blow-up with size $s$ with respect to $\mathcal{P}$.*

**Proof.** Set $a := \lceil 2es^2D \rceil$. For each $i \in [r]$, we repeatedly prune vertices of largest degree in $V_i$, stopping when exactly $a$ vertices remain. We refer to these reduced sets as $V_i'$. It is straightforward to see that $V_i'$ is a subset of $V_i$ for all $i$ in $[r]$. Now, consider the graph $G'$ induced by the union of these reduced sets: $G' = G[V_1' \cup V_2' \cup \ldots \cup V_r']$. Notably, any independent transversal blow-up in $G'$ is also valid in $G$. Hence, we can assume without loss of generality that each $V_i$ contains exactly $a$ vertices. Next, we proceed to randomly select $s$ vertices from each $V_i$, doing so independently and uniformly. Let $W$ denote the resulting random set of vertices. To establish the proof, we aim to demonstrate that $W$ forms an independent set in $G$ with positive probability.

For any edge $f$ in $G$, define the event $A_f$ as the occurrence where both endpoints of $f$ are included in $W$. It is immediate that

$$\mathbb{P}[A_f] \le \left( \frac{\binom{a-1}{s-1}}{\binom{a}{s}} \right)^2 = \frac{s^2}{a^2}.$$

Furthermore, if $f$ connects vertices in $V_i$ and $V_j$, then $A_f$ is independent of any event associated with edges whose endpoints do not lie within $V_i \cup V_j$. This gives rise to a dependency digraph where the maximum degree is bounded by $2aD - 2$. Given that $e \cdot 2aD \cdot \frac{s^2}{a^2} = \frac{2es^2D}{a} < 1$, we can invoke Lemma 2.3 to conclude that, with positive probability, none of the events $A_f$ happens. Consequently, $W$ constitutes an independent set that includes $s$ vertices from each $V_i$, thereby concluding the proof. □

**Lemma 4.2.** [16, Kang and Kelly] *Given a multipartite graph $G$ with a partition $\mathcal{P} = (V_1, V_2, \ldots, V_r)$, if $\overline{D}_{\mathcal{P}}(G) \leq D$ and $\mathcal{P}$ is $4D$-thick, then $G$ has an independent transversal with respect to $\mathcal{P}$.*

**Remark 1.** *The above result surpasses the one provided by Theorem 4.1 for $s = 1$ and is currently the best known. However, the tightness of the constant 4 remains unknown.*

**Theorem 4.3.** *Given a multipartite graph $G$ with a partition $\mathcal{P} = (V_1, V_2, \ldots, V_r)$, if $\overline{D}_{\mathcal{P}}(G) \leq D$ and $\mathcal{P}$ is $4s^2D$-thick, then $G$ has an independent transversal blow-up with size $s$ with respect to $\mathcal{P}$.*

**Proof.** Take the auxiliary graph $\dot{G}$ based on $G$ and $\mathcal{P}$. According to Definition 1, $\dot{G}$ is a multipartite graph with a partition $\dot{\mathcal{P}} = (\dot{V}_1, \dot{V}_2, \ldots, \dot{V}_r)$. By Observation 1, we know that $\overline{D}_{\dot{\mathcal{P}}}(\dot{G}) \leq s\binom{4s^2D-1}{s-1}D$ and $\dot{\mathcal{P}}$ is $\binom{4s^2D}{s}$-thick. Since

$$\binom{4s^2D}{s} = \binom{4s^2D-1}{s-1}\frac{4s^2D}{s} = 4s\binom{4s^2D-1}{s-1}D,$$

we can utilize Lemma 4.2 to deduce that $\dot{G}$ possesses an independent transversal with respect to $\dot{\mathcal{P}}$. By referring once again to Definition 1, this independent transversal can be promptly converted into an independent transversal blow-up of size $s$ in $G$ with respect to $\mathcal{P}$. □

Now we shift our attention to the maximum degree setting. Below is Haxell's theorem, which we previously mentioned in the introduction.

**Lemma 4.4.** [13, Haxell] *Given a multipartite graph $G$ with a partition $\mathcal{P} = (V_1, V_2, \ldots, V_r)$, if $\Delta(G) \leq \Delta$ and $\mathcal{P}$ is $2\Delta$-thick, then $G$ has an independent transversal with respect to $\mathcal{P}$.*

Using this Haxell's theorem, we obtain a comparable result in the context of maximum degree. Specifically, Haxell's theorem provides a condition under which a graph-partition instance $(G, \mathcal{P})$ with bounded maximum degree admits an independent transversal. In the context of our discussion, this implies that if the maximum degree in the auxiliary graph $\dot{G}$ (constructed according to Definition 1) is bounded, then we can find an independent transversal in $\dot{G}$ with respect to $\dot{\mathcal{P}}$. This, in turn, can be translated back to an independent transversal blow-up of size $s$ in $G$ with respect to $\mathcal{P}$, similar to the previous result based on Lemma 4.2. Thus, Haxell's theorem allows us to draw a parallel conclusion in the setting of maximum degree. We omit the proof and leave it for interested readers to explore.

**Theorem 4.5.** *Given a multipartite graph $G$ with a partition $\mathcal{P} = (V_1, V_2, \ldots, V_r)$, if $\Delta(G) \le \Delta$ and $\mathcal{P}$ is $2s^2\Delta$-thick, then $G$ has an independent transversal blow-up with size $s$ with respect to $\mathcal{P}$.*

As mentioned earlier, Glock and Sudakov [10], and independently, Kang and Kelly [16] derived the following theorem from the perspective of a local degree and an average block degree.

**Lemma 4.6.** [10, Glock and Sudakov], [16, Kang and Kelly] *For every $\varepsilon > 0$, there exists $\gamma > 0$ such that the following hold. If $G$ is a multipartite graph with a partition $\mathcal{P} = (V_1, V_2, \ldots, V_r)$, where $\overline{D}_{\mathcal{P}}(G) \le D$, $\Gamma_{\mathcal{P}}(G) \le \gamma D$, and $\mathcal{P}$ is $(1+\varepsilon)D$-thick, then $G$ has an independent transversal with respect to $\mathcal{P}$.*

To extend their results to the independent transversal blow-up, we present an intuitive result with a quite simple proof.

**Theorem 4.7.** *For every $\varepsilon > 0$, there exists $\gamma > 0$ such that the following hold. If $G$ is a multipartite graph with a partition $\mathcal{P} = (V_1, V_2, \ldots, V_r)$, where $\overline{D}_{\mathcal{P}}(G) \le D$, $\Gamma_{\mathcal{P}}(G) \le \gamma D$, and $\mathcal{P}$ is $(s^2+\varepsilon)D$-thick, then $G$ has an independent transversal blow-up with size $s$ with respect to $\mathcal{P}$.*

**Proof.** The proof bears a resemblance to that of Theorem 4.3. Define $\dot{D}$ as $s\binom{(s^2+\varepsilon)D-1}{s-1}D$ and denote by $\gamma$ the value derived from Lemma 4.6, where $\frac{\varepsilon}{s^2}$ is substituted for $\varepsilon$.

By examining the auxiliary graph $\dot{G}$ and applying Observation 1, we conclude that $\overline{D}_{\dot{\mathcal{P}}}(\dot{G}) \le \dot{D}$ and $\Gamma_{\dot{\mathcal{P}}}(\dot{G}) \le \gamma \dot{D}$. Furthermore, each block in $\dot{G}$ has a size of at least:

$$\binom{(s^2+\varepsilon)D}{s} = \binom{(s^2+\varepsilon)D-1}{s-1}\frac{s^2+\varepsilon}{s}D = \left(1+\frac{\varepsilon}{s^2}\right)\dot{D}.$$

Consequently, according to Lemma 4.6, $\dot{G}$ possesses an independent transversal with respect to $\dot{\mathcal{P}}$. This, in accordance with Definition 1, can be translated into an independent transversal blow-up of size $s$ in $G$ with respect to $\mathcal{P}$. □

## 4.2 Sharpen the bound

In this section, our objective is to refine Theorem 4.7 by improving the partition thickness from $(s^2+\varepsilon)D$ to $(s+\varepsilon)D$, as elaborated in Theorem 4.8. Furthermore, in Theorem 4.17, we demonstrate the tightness of the number $(s+\varepsilon)D$ by proving that it cannot be reduced to $sD$ without losing the validity of the theorem.

**Theorem 4.8.** *For every $\varepsilon > 0$, there exists $\gamma > 0$ such that the following hold. If $G$ is a multipartite graph with a partition $\mathcal{P} = (V_1, V_2, \ldots, V_r)$, where $\overline{D}_{\mathcal{P}}(G) \le D$, $\Gamma_{\mathcal{P}}(G) \le \gamma D$, and $\mathcal{P}$ is $(s+\varepsilon)D$-thick, then $G$ has an independent transversal blow-up with size $s$ with respect to $\mathcal{P}$.*

**Lemma 4.9.** [16, Kang and Kelly] *For every $\varepsilon > 0$, there exists $\gamma_0, D_0 > 0$ such that the following holds for all $\gamma < \gamma_0$ and $D > D_0$. For every multipartite graph $G$ with a partition $\mathcal{P} = (V_1, V_2, \ldots, V_r)$ satisfying*

- *$|V_i| \ge (s+\varepsilon)D$;*
- *$\overline{D}_{\mathcal{P}}(G) \le D$;*
- *$\Gamma_{\mathcal{P}}(G) \le \gamma D$,*

*there exists $D' \geq D$ and subsets $W_i \subset V_i$ for each $1 \leq i \leq r$ such that the $r$-partite subgraph $G'$ of $G$ induced by the set $\cup W_i$ satisfying the following properties:*

- $|W_i| \geq (s + \varepsilon/2)D'$;
- $\overline{D}_{\mathcal{P}}(G') \leq D'$;
- $\Gamma_{\mathcal{P}}(G') \leq D'^{1/5}$.

**Remark 2.** *Kang and Kelly validated this lemma for $s = 1$ (refer to [16, Theorem 4]). However, their proof can be adapted to hold for all values of $s$, and thus, we refrain from elaborating on the specifics here.*

**Lemma 4.10.** *For every $\varepsilon > 0$, the following holds . For every multipartite graph $G$ with a partition $\mathcal{P} = (V_1, V_2, \ldots, V_r)$, where $\overline{D}_{\mathcal{P}}(G) \leq D$ and $\mathcal{P}$ is $(s + \varepsilon/2)D$-thick, there exist $D' \geq D$ and subsets $W_i \subset V_i$ for each $1 \leq i \leq r$ such that the $r$-partite subgraph $G'$ of $G$ induced by the set $\bigcup W_i$ satisfying the following properties:*

- $|W_i| \geq (s + \varepsilon/4)D'$;
- $\overline{D}_{\mathcal{P}'}(G') \leq D'$;
- $\Delta(G') \leq 16sD'/\varepsilon$.

**Proof.** For each $i \in [r]$, define

$$B_i = \left\{ v \in V_i : \ d_G(v) \geq \frac{16sD}{\varepsilon} \right\}.$$

It follows

$$|B_i| \cdot \frac{16sD}{\varepsilon} \leq \sum_{v \in V_i} d_G(v) \leq D|V_i|,$$

from which we can deduce that $|B_i| \leq \frac{\varepsilon |V_i|}{16s}$. Let $W_i := V_i \setminus B_i$ ($i \in [r]$) and let $G'$ be the $r$-partite subgraph of $G$ induced by the set $\bigcup W_i$. Set $D' = \frac{D}{1 - \frac{\varepsilon}{16s}}$, we have

- $|W_i| \geq \left(1 - \frac{\varepsilon}{16s}\right)|V_i| \geq \left(1 - \frac{\varepsilon}{16s}\right)\left(s + \frac{\varepsilon}{2}\right) D \geq \left(s + \frac{\varepsilon}{4}\right) D'$;
- $\overline{d}_{G'}(W_i) \leq \frac{|V_i|}{|W_i|}\overline{d}_G(V_i) \leq \frac{D}{1 - \frac{\varepsilon}{16s}} = D'$;
- $\Delta(G') \leq \frac{16sD}{\varepsilon} = \frac{16s}{\varepsilon}\left(1 - \frac{\varepsilon}{16s}\right) D' \leq \frac{16sD'}{\varepsilon}$,

as desired. □

Proving the following assertion, with the support of Lemmas 4.9 and 4.10, will suffice for proving Theorem 4.8.

**Theorem 4.11.** *For every $\varepsilon > 0$, there exists $D_0 > 0$ such that the following holds for all $D > D_0$. For every multipartite graph $G$ with a partition $\mathcal{P} = (V_1, V_2, \ldots, V_r)$, if $\overline{D}_{\mathcal{P}}(G) \leq D$, $\Delta(G) \leq 16sD/\varepsilon$, $\Gamma_{\mathcal{P}}(G) \leq D^{1/5}$, and $\mathcal{P}$ is $(s + \varepsilon/4)D$-thick, then $G$ has an independent transversal blow-up with size $s$ with respect to $\mathcal{P}$.*

We find a partial independent transversal blow-up with size $s$ (PITS for short) as follows.

In the first round, we start with $G := G(1)$ and the partition $\mathcal{P} := \mathcal{P}(1)$, accompanied by an empty PITS. Subsequently, we technically identify a PITS(1) with respect to $\mathcal{P}$ technically. Next, we outline the framework of this technical process.

Now, let us assume we are at the start of the $t$-th round, possessing a multipartite graph $G(t)$ with a $S(t)$-thick partition $\mathcal{P}(t) = (V_{i_1}(t), V_{i_2}(t), \ldots, V_{i_{r(t)}}(t))$ and with maximum block average degree at most $D(t)$.

Set $p := \frac{1}{\log^3 D}$. We define the sequences $\{S(t)\}$ and $\{D(t)\}$ recursively as follows:

$$S(1) = (s+\varepsilon)D,\ S(t+1) = \left(1 - p + \frac{\varepsilon p}{2s}\right)S(t),$$
$$D(1) = D,\ D(t+1) = \left(1 - p + \frac{\varepsilon p}{4s}\right)D(t),$$

**Observation 3.** *Set $t^* := \frac{8s\log(2es)}{\varepsilon p}$. For all $1 \le t \le t^*$, the following hold.*

1. *The sequence $\left\{\frac{D(t+1)}{S(t+1)}\right\}$ is decreasing with respect to $t$.*
2. *$S(t+1) = \Theta(D)$ and $D(t+1) = \Theta(D)$.*
3. *$S(t^*+1) \ge 2es^2 D(t^*+1)$*

**Proof.** 1. It is immediate that

$$\frac{D(t+1)}{S(t+1)} = \left(\frac{1 - p + \frac{\varepsilon p}{4s}}{1 - p + \frac{\varepsilon p}{2s}}\right)\frac{D(t)}{S(t)} \le \frac{D(t)}{S(t)} \le \cdots \le \frac{D(1)}{S(1)}.$$

2. The conclusion $S(t+1) = \Theta(D)$ follows due to the fact that

$$\begin{aligned} S(t+1) &= \left(1 - p + \frac{\varepsilon p}{2s}\right)^t S(1) < S(1) = O(D), \\ S(t+1) &= \left(1 - p + \frac{\varepsilon p}{2s}\right)^t S(1) \ge \exp\left(-(2-\varepsilon/s)pt\right) S(1) \ge \exp\left(-(2-\varepsilon/s)pt^*\right) S(1) \\ &\ge s\exp\left(-\frac{8(2-\varepsilon/s)\log(2es)}{\varepsilon}\right) D = \Omega(D). \end{aligned}$$

Similarly, we have $D(t+1) = \Theta(D)$. Here we apply the inequality $1 - x \ge e^{-2x}$ for all $0 \le x \le 0.5$.

3. A simple calculation implies

$$\begin{aligned} \frac{D(t^*+1)}{S(t^*+1)} &= \left(\frac{1 - p + \frac{\varepsilon p}{4s}}{1 - p + \frac{\varepsilon p}{2s}}\right)^{t^*}\frac{D(1)}{S(1)} < \left(1 - \frac{\frac{\varepsilon p}{4s}}{1 - p + \frac{\varepsilon p}{2s}}\right)^{t^*}\frac{1}{s} \\ &\le \frac{1}{s}\exp\left(-\frac{\frac{\varepsilon p t^*}{4s}}{1 - p + \frac{\varepsilon p}{2s}}\right) \le \frac{1}{s}\exp\left(-\frac{2}{1 - p + \frac{\varepsilon p}{2s}}\log(2es)\right) \\ &\le \frac{1}{s}\exp\left(-\log(2es)\right) = \frac{1}{2es^2}, \end{aligned}$$

as desired. □

During the $t$-th round, we aim to do the following:

**A0** Continuously remove vertices possessing the highest degree from each $V_i(t)$ until precisely $S(t)$ vertices remain, resulting in an updated graph-partition pair, still denoted as $(G(t), \mathcal{P}(t))$, where $\mathcal{P}(t) = \bigcup_{i \in I(t)} V_i(t)$. Here, $I(t)$ represents the set of indices corresponding to the partition $\mathcal{P}(t)$, in other words, each index $i \in I(t)$ uniquely identifying a subset $V_i(t)$ within the partition $\mathcal{P}(t)$;

**A1** find a $s$-set independently uniformly at random from each block of the partition $\mathcal{P}(t)$, and collect those $s$-sets into a set $T$;

**A2** from $T$, find a subset $T'$ such that $\cup T'$ forms an independent set of $G(t)$, and gather the indices of the blocks from which each set in $T'$ originates into a set $J'$;

**A3** eliminate the entire block $V_j(t)$ with $j \in J'$ and remove all vertices that are neighbors of any vertex in $\cup T$,

ensuring that we arrive at a multipartite graph $G(t+1)$ with a partition $\mathcal{P}(t+1) = \bigcup_{i \in I(t+1)} V_i(t+1)$ such that for each $i \in I(t+1)$,

**G1** $|V_i(t+1)| \geq S(t+1)$;

**G2** $\overline{d}_{G(t+1)}(V_i(t+1)) \leq D(t+1)$.

Note that the operation outlined in **A0** guarantees that the average degree within each block and the local degree of the graph with respect to the partition do not increase. Additionally, it confirms that the partition retains an $S(t)$-thick structure, which will significantly simplify subsequent calculations. Following this, the set $T'$ obtained during step **A2** forms a PITS, denoted as PITS($t$). The operation described in **A3** guarantees that the combined sets of PITS(1), ..., PITS($t$) constitute a PITS of $G$. Subsequently, we proceed to the $(t+1)$-th round, starting with a $S(t+1)$-thick partition $\mathcal{P}(t+1)$ of a multipartite graph $G(t+1)$, where the maximum block average degree at most $D(t+1)$. This allows us to proceed with the iterative process outlined above, continuing until we reach the $(t^*+1)$-th step. Since $S(t^*+1) \geq 2es^2 D(t^*+1)$ by Observation 3(3), applying Theorem 4.1 to the graph-partition pair $(G(t^*+1), \mathcal{P}(t^*+1))$, we obtain an independent transversal blow-up with size $s$ (shortly ITS) in $G(t^*+1)$ with respect to $\mathcal{P}(t^*+1)$. Combining this ITS with the union of PITS(1), ..., PITS($t^*$) that we already possessed, we ultimately obtain an ITS for the entire graph $G$.

From this point forward, we shall delve into more specific details concerning **A1** and **A2**, as the achievement of **G1** and **G2** entirely depends on how **A1** and **A2** are specifically designed.

For **A1**, our goal is to identify an $s$-set within each block $V_i(t)$ of the partition $\mathcal{P}(t)$ that possesses a "good" property. To expedite this task, we employ the auxiliary graph $\dot{G}(t)$, which is derived from $G(t)$ and defined in Definition 1. We reformulate Observation 1 as follows.

**Observation 4.** *For $\dot{G}(t)$ and $\dot{\mathcal{P}}(t) = \bigcup_i \dot{V}_i(t)$, the following holds.*

1. $|\dot{V}_i(t)| = \binom{S(t)}{s}$;
2. $\Delta(\dot{G}(t)) \le s\binom{S(t)-1}{s-1}\Delta(G(t))$;
3. $\Gamma_{\dot{\mathcal{P}}(t)}(\dot{G}(t)) \le s\binom{S(t)-1}{s-1}\Gamma_{\mathcal{P}(t)}(G(t))$;
4. $\overline{D}_{\dot{\mathcal{P}}(t)}(\dot{G}(t)) \le s\binom{S(t)-1}{s-1}\overline{D}_{\mathcal{P}(t)}(G(t))$.

We independently and randomly select a vertex $v_{Z_i}$ from each $\dot{V}_i(t)$, and subsequently construct a random set $\dot{T}$ by taking the union of all these vertices: $\dot{T} = \bigcup_i \{v_{Z_i}\}$.

**Claim 4.12.** *With a probability of at least $1 - \frac{1}{2}\exp(-\log D \log\log D)$, each vertex of $\dot{G}(t)$ will have no more than $\log^2 D$ neighboring vertices within $\dot{T}$.*

**Proof.** Fix a vertex $v \in V(\dot{G}(t))$. For each $i$, let $X_i$ be the indicator that has value 1 if $v_{s_i}$ is a neighbor of $v$ in $\dot{G}(t)$, and 0 otherwise. Clearly, $\mathbb{P}[X_i = 1] = \frac{|N_{\dot{G}(t)}(v)\cap \dot{V}_i(t)|}{|\dot{V}_i(t)|}$.

Then, $X := \sum_i X_i$ is the number of neighbors of $v$ in $\dot{T}$, a sum of independent Bernoulli random variables. It follows

$$\mathbb{E}[X] = \sum_i \mathbb{P}[X_i = 1] = \sum_i \frac{|N_{\dot{G}(t)}(v) \cap \dot{V}_i(t)|}{|\dot{V}_i(t)|}$$
$$\overset{\text{Obs. 4(1)}}{\le} \frac{d_{\dot{G}(t)}(v)}{\binom{S(t)}{s}} \overset{\text{Obs. 4(2)}}{\le} \frac{s\binom{S(t)-1}{s-1}\Delta(G(t))}{\binom{S(t)}{s}} \le O\left(\frac{\Delta(G(t))}{S(t)}\right) = O(1).$$

Now, applying Chernoff's Bound (Lemma 2.1), we have

$$\mathbb{P}[X \ge \log^2 D] \le \exp(-\log^2 D) \le \frac{1}{2}\exp(-\log D \log\log D),$$

as desired. □

Based on Claim 4.12, we had earlier pinpointed a specific set $\dot{T} = \bigcup_i \{v_{Z_i}\}$. Considering the property of $\dot{T}$ described in Claim 4.12, we then proceed to identify a subset $\dot{T}' \subset \dot{T}$ that gives rise to an independent set within $\dot{G}(t)$. To accomplish this, we transition our randomized algorithm into its second phase:

**A2-1** for each block $\dot{V}_i$ of $\dot{\mathcal{P}}(t)$, activate it with probability $p := \frac{1}{\log^3 D}$; we collect the indices of all activated blocks $\dot{V}_i$ into $J$;

**A2-2** identify a maximal independent set $\dot{T}'$ from $\bigcup_{i\in J}\{v_{Z_i}\}$; then, define the set $J'$ as containing all $i$ such that $v_{Z_i}$ belongs to $\dot{T}'$.

Now, for **A2**, we have $T' := \bigcup \dot{T}'$ and the desired set $J'$.

Next, we evaluate how **A3** would alter the structure of $G(t)$, which properties would be inherited into the next round, and thus how **G1** and **G2** would be achieved.

For $i \in J'$, the block $V_i(t)$ is completely removed from $G(t)$, prompting us to examine the alteration in the size of blocks $V_i(t)$ where $i \notin J'$. Let $v$ be any vertex in $V_i(t)$. According to **A3**, $v$ is deleted only if $v$ is adjacent to some vertex in $T'$. For each neighbor $u$ of $v$, the probability that $u$ is included in $T'$ is at most $p \cdot s/S(t)$. By a simple union bound, we have

$$\begin{aligned}\mathbb{P}[v \in V_i(t+1)] = \mathbb{P}[N_{G(t)}(v) \cap T' = \emptyset] &= 1 - \mathbb{P}[N_{G(t)}(v) \cap T' \neq \emptyset] \\ &\geq 1 - \frac{d_{G(t)}(v)ps}{S(t)} := p_v.\end{aligned}$$

In other words, a vertex $v \in V_i(t)$ will be inherited into the $(t+1)$-th round with probability at least $p_v$. However, for the convenience of the following analysis, we may introduce artificial intervention on the basis of **A3** to determine whether to delete $v$ so that the probability of it being retained is exactly $p_v$. Keeping this in mind, we introduce an artificial Bernoulli random variable $B_v$, which is independent of all other variables, such that for each $i \notin J'$ and each $v \in V_i(t)$,

$$\mathbb{P}[N_{G(t)}(v) \cap T' = \emptyset \wedge B_v = 1] = p_v.$$

This serves our purpose adequately. In this context, $B_v = 0$ signifies that we have manually deleted $v$ as an artificial step.

Therefore, for any index $i \in I(t) \setminus J'$ and subset $U \subseteq V_i(t)$, we have

$$\mathbb{E}[|U \cap V_i(t+1)|] = \sum_{v \in U} p_v.$$

In the following, we show that $|U \cap V_i(t+1)|$ is concentrated around its mean.

**Claim 4.13.** *For any index $i \in I(t) \setminus J'$ and subset $U \subseteq V_i(t)$, let $U' = U \cap V_i(t+1)$ and $Y = |U'|$. We have*

$$\mathbb{P}\left[\left|Y - \sum_{v \in U} p_v\right| \geq D^{2/5} \log^4 D \sqrt{p|U|}\right] \leq \exp(-\log D \log\log D).$$

**Proof.** Rather than proving the concentration of $Y$ directly, we demonstrate that $R := |U \setminus V_i(t+1)|$ is concentrated. Here $R$ counts the number of vertices deleted from $V_i(t)$ and is completely dependent on a series of independent trials that indicate whether a block is activated or whether a vertex is deleted artificially. Clearly,

$$\mathbb{E}[R] = \sum_{v \in U}(1 - p_v) = \sum_{v \in U} \frac{d_{G(t)}(v)ps}{S(t)} \leq \sum_{v \in U} \frac{O(D)ps}{S(t)} = O(p|U|).$$

For any vertex $v$ that is removed from $V_i(t)$, there would be an activated block, say $V_j(t)$, such that the selected $s$-set $Z_j$ in that block has one vertex neighboring to $v$, or $B_v = 0$. Since these events certify the removal of $v$, the random variable $R$ is 1-certifiable. On the other hand, inverting the activation status of a block or altering the value of a Bernoulli random variable $B_v$ would affect $R$ by at most $sD^{1/5}$, constrained by the bound of the local degree. This implies that the random variable $R$ is $sD^{1/5}$-Lipschitz.

Set $t = D^{2/5}\log^4 D\sqrt{p|U|}/2$. Since $t \geq 20sD^{1/5}\sqrt{\mathbb{E}[R]} + 64s^2D^{2/5}$, we can apply Talagrand's Inequality (Lemma 2.2) with parameter $t$, obtaining

$$\mathbb{P}\left[|R - \mathbb{E}[R]| \geq D^{2/5}\log^4 D\sqrt{p|U|}\right] \leq 4\exp\left(-\frac{t^2}{8s^2D^{2/5}(\mathbb{E}[R]+t)}\right)$$
$$\leq O\left(\exp\left(-\frac{D^{4/5}\log^8 Dp|U|}{D^{2/5}p|U| + D^{4/5}\log^4 D\sqrt{p|U|}}\right)\right) \leq O\left(\exp\left(-\frac{D^{2/5}\log^8 D\sqrt{p|U|}}{\sqrt{p|U|} + D^{2/5}\log^4 D}\right)\right)$$
$$\leq O\left(\exp\left(-\frac{D^{2/5}\log^8 D\sqrt{p}}{\sqrt{p} + D^{2/5}\log^4 D}\right)\right) \leq O\left(\exp\left(-\frac{D^{2/5}\log^8 D}{1 + D^{2/5}\log^{5.5} D}\right)\right) \leq O\left(\exp\left(-\frac{D^{2/5}\log^8 D}{2D^{2/5}\log^{5.5} D}\right)\right)$$
$$\leq O\left(\exp(-\log^{2.5} D)\right) \leq O\left(\exp(-\log^2 D)\right) \leq \exp(-\log D\log\log D).$$

Since $Y + R = |U| = \mathbb{E}[R] + \sum_{v\in U} p_v$, the claim holds now based on the inequality derived above. □

**Claim 4.14.** *For every vertex $v \in V(G(t))$ with $d_{G(t)}(v) \geq D^{2/5}\log^{12} D$, with a probability of at least $1 - \exp(-\log D\log\log D)$, we have*

$$d_{G(t+1)}(v) \leq \left(1 - p + \frac{\varepsilon p}{5s}\right)d_{G(t)}(v).$$

**Proof.** Since $v$ can be inherited into the $(t+1)$-round (otherwise it is not in our case), $v \in V_i(t)$ for some $i \notin J'$. For each $i \in I(t)$, let $c_i$ be the number of neighbors of $v$ in $V_i(t)$ and let $K := \{i \mid c_i > 0\}$. Since $V_i(t)$ with $j \in J'$ will be completely removed as a result of **A3**, the vertex $v$ will lose a degree of at least $\sum_{i\in K} \mathbb{1}_{i\in J'}c_i$ from $G(t)$ to $G(t+1)$, i.e.,

$$d_{G(t)}(v) - d_{G(t+1)}(v) \geq \sum_{i\in K} \mathbb{1}_{i\in J'}c_i.$$

Define the random variable

$$Z := d_{G(t)}(v) - \sum_{i\in K} \mathbb{1}_{i\in J'}c_i = \sum_{i\in K} \mathbb{1}_{i\notin J'}c_i.$$

It is sufficient to show

$$\mathbb{P}\left[Z \leq \left(1 - p + \frac{\varepsilon p}{5s}\right)d_{G(t)}(v)\right] \geq 1 - \exp(-\log D\log\log D).$$

Using Claim 4.12, we first infer that the following event $\Psi$ hold with probability

$$\mathbb{P}[\Psi] \geq 1 - \frac{1}{2}\exp(-\log D\log\log D):$$

- *each vertex of $\dot{G}(t)$ has at most $\log^2 D$ neighboring vertices within $\dot{T} = \bigcup_i\{v_{Z_i}\}$, where $v_{Z_i}$ is an independently and randomly selected vertex from $\dot{V}_i(t)$.*

Suppose $\Psi$ holds and fix any choice $\dot{T}$ for which $\Psi$ holds. In order to show that $Z$ is unlikely to be too large when $\Psi$ holds, we shall bound the conditional probability $\mathbb{P}[i \notin J' \mid \dot{T}]$.

For each $i$, let $Q_i$ be the set of indices $j$ such that $v_{Z_j}$ is adjacent to $v_{Z_i}$ in $\dot{G}(t)$. By Claim 4.12, $|Q_i| \leq \log^2 D$. According to **A2.1** and **A2.2**, $i \in J'$ if and only $\dot{V}_i$ is activated and $\dot{V}_j$ with $j \in Q_i$ is not activated, which happens with probability at least $p(1-p)^{|Q_i|} \geq p(1-p|Q_i|) \geq p(1-p\log^2 D)$. Therefore,

$$\mathbb{P}[i \notin J' \mid \dot{T}] = 1 - \mathbb{P}[i \in J' \mid \dot{T}] \leq 1 - p + p^2 \log^2 D$$

and

$$\mathbb{E}[Z \mid \dot{T}] \leq (1 - p + p^2 \log^2 D) d_{G(t)}(v).$$

Define $\bar{K}$ as the collection of indices $i$ satisfying the condition that either $i \in K$ or there exists an index $j$ such that the number of edges between $\dot{V}_i(t)$ and $\dot{V}_j(t)$ is positive. Now, the random variable $Z$ is only determined by the activation status of the parts $\dot{V}_i(t)$ for indices $i \in \bar{K}$. Furthermore, the activation of a particular $\dot{V}_i(t)$ influences the event $\{j \notin J'\}$ only if $j \in Q_i \cup \{i\}$. Therefore, toggling the activation state of any block can alter $Z$ by an amount bounded by $D^{1/5}(1+|Q_i|) \leq 2D^{1/5}\log^2 D$. This bound arises because, according to the local degree condition, each vertex $v$ has at most $D^{1/5}$ neighboring vertices within any block of $G(t)$. Consequently, $Z$ satisfies a Lipschitz condition with a constant of $2D^{1/5}\log^2 D$. Ultimately, the 1-certifiability of $Z$ becomes evident, as the occurrence of the event $\{i \notin J'\}$ can be verified through either the non-activation of the block $\dot{V}_i(t)$ or the activation of a block $\dot{V}_j(t)$ where $j \in Q_i$.

Set $t = \frac{\varepsilon p d_{G(t)}(v)}{10s}$. First, $t \geq 40D^{1/5}\log^2 D\sqrt{\mathbb{E}[Z \mid \dot{T}]} + 256D^{2/5}\log^4 D$. This holds if

$$\varepsilon d_{G(t)}(v) \geq 400sD^{1/5}\log^5 D\sqrt{d_{G(t)}(v)} + 2560sD^{2/5}\log^7 D.$$

Next, we apply Talagrand's inequality (Lemma 2.2) with this parameter $t$, concluding

$$\mathbb{P}\left[Z > \left(1 - p + \frac{\varepsilon p}{5s}\right) d_{G(t)}(v) \,\middle|\, \dot{T}\right] \leq \mathbb{P}\left[Z > \frac{\varepsilon p d_{G(t)}(v)}{5s} \,\middle|\, \dot{T}\right]$$
$$\leq 4\exp\left(-\frac{t^2}{32D^{2/5}\log^4 D(\mathbb{E}[Z \mid \dot{T}] + t)}\right) \leq O\left(\exp\left(-\frac{p^2 d_{G(t)}(v)}{D^{2/5}\log^4 D(1 - p + p^2\log^2 D) + pD^{2/5}\log^4 D}\right)\right)$$
$$\leq O\left(\exp\left(-\frac{d_{G(t)}(v)}{D^{2/5}\log^{10} D}\right)\right) \leq O\left(\exp\left(-\log^2 D\right)\right) \leq \frac{1}{2}\exp\left(-\log D \log\log D\right)$$

provided $d_{G(t)}(v) \geq D^{2/5}\log^{12} D$.

Since this bound holds for any choice $\dot{T}$ satisfying $\Psi$, we conclude that

$$\mathbb{P}\left[Z \leq \left(1 - p + \frac{\varepsilon p}{5s}\right) d_{G(t)}(v) \,\middle|\, \Psi\right] \geq 1 - \frac{1}{2}\exp\left(-\log D \log\log D\right).$$

Therefore,

$$\mathbb{P}\left[Z \leq \left(1 - p + \frac{\varepsilon p}{5s}\right) d_{G(t)}(v)\right] = \mathbb{P}[\Psi] \cdot \mathbb{P}\left[Z \leq \left(1 - p + \frac{\varepsilon p}{5s}\right) d_{G(t)}(v) \,\middle|\, \Psi\right]$$
$$\geq \left(1 - \frac{1}{2}\exp\left(-\log D \log\log D\right)\right)^2 \geq 1 - \exp\left(-\log D \log\log D\right),$$

as desired. □

We will utilize the Lovász Local Lemma (Lemma 2.3) to prove that both **G1** and **G2** hold with a positive probability. Initially, we have

$$\mathbb{P}\left[|V_i(t+1)| < S(t+1)\right] = \mathbb{P}\left[|V_i(t+1)| < \left(1 - p + \frac{\varepsilon p}{2s}\right)S(t)\right] \le \mathbb{P}\left[|V_i(t+1)| < \left(1 - p + \frac{\varepsilon p}{2s}\right)|V_i(t)|\right]$$
$$= \mathbb{P}\left[|V_i(t+1)| < \left(1 - \frac{ps}{s+\varepsilon}\right)|V_i(t)| - \left(\varepsilon p\sqrt{|V_i(t)|}\left(\frac{1}{s+\varepsilon} - \frac{1}{2s}\right)\right)\sqrt{|V_i(t)|}\right].$$

Since

$$\sum_{v\in V_i(t)} p_v = \sum_{v\in V_i(t)}\left(1 - \frac{d_{G(t)}(v)ps}{S(t)}\right) = |V_i(t)| - \frac{ps}{S(t)}\sum_{v\in V_i(t)} d_{G(t)}(v) = |V_i(t)|\left(1 - \frac{ps}{S(t)}\bar{d}_{G(t)}(V_i(t))\right)$$
$$\ge |V_i(t)|\left(1 - ps\frac{D(t)}{S(t)}\right) \overset{\text{Obs. 3(1)}}{\ge} |V_i(t)|\left(1 - ps\frac{D(1)}{S(1)}\right) = |V_i(t)|\left(1 - \frac{ps}{s+\varepsilon}\right),$$

and

$$\varepsilon p\sqrt{|V_i(t)|}\left(\frac{1}{s+\varepsilon} - \frac{1}{2s}\right) \ge \varepsilon p\sqrt{S(t)}\left(\frac{1}{s+\varepsilon} - \frac{1}{2s}\right) \overset{\text{Obs. 3(2)}}{\ge} \Theta(pD^{1/2}) \ge D^{2/5}\log^3 D,$$

continuing from the previous calculation, we further obtain

$$\mathbb{P}\left[|V_i(t+1)| < S(t+1)\right] \le \mathbb{P}\left[|V_i(t+1)| < \sum_{v\in V_i(t)} p_v - D^{2/5}\log^3 D\sqrt{|V_i(t)|}\right]$$
$$\le \exp(-\log D \log\log D)$$

by Claim 4.13, where we substitute $U = V_i(t)$. Therefore, we conclude for each $i \in I(t+1)$ that

**F1** *the bad event $\mathcal{B}_i$ that* **G1** *does not hold for $i$ occurs with probability at most*

$$\exp(-\log D \log\log D).$$

We now turn to check **G2** conditioning on **G1**. Set $\mu = \lceil D^{9/10}\rceil$. For each $i \in I(t)$, we categorize the vertices of $V_i(t)$ based on their degrees. We define:

$$X_j := \{v \in V_i(t) \mid j\mu \le d_{G(t)}(v) < (j+1)\mu\}, \quad X'_j := X_j \cap V_i(t+1).$$

Since $d_{G(t)}(v) \le \Delta(G(t)) = O(D)$, $j = O(D/\mu) \le O(D^{1/10})$.

By substituting $U = X_j$ in Claim 4.13 and applying Claim 4.14, we deduce, using the union bound, that with probability at least $1 - O(D)\exp(-\log D\log\log D)$, the following statements hold:

**(a)** $\left|X'_j\right| = \sum_{v\in X_j} p_v \pm D^{2/5}\log^3 D\sqrt{\left|X_j\right|}$;

**(b)** $d_{G(t+1)}(v) \le \left(1 - p + \frac{\varepsilon p}{5s}\right)d_{G(t)}(v)$ for every vertex $v \in V_i(t) \setminus X_0$.

Our objective is to establish that for each $i \in I(t+1)$, the following inequality holds:

$$\frac{\sum_{v \in V_i(t+1)} d_{G(t)}(v)}{|V_i(t+1)|} \le (1 + D^{-1/300}) D(t). \tag{$\spadesuit$}$$

Prior to proving this inequality, we will first illustrate how it serves our purpose.

Utilizing the following facts:

$$|X_0'| \le |V_i(t)| \overset{\mathbf{A1}}{=} |S(t)| \le O(D),$$
$$|V_i(t+1)| D(t) \overset{\mathbf{G1}}{\ge} S(t+1) D(t) \overset{\text{Obs. 3(2)}}{\ge} \Theta(D^2),$$

along with the statement **(b)**, we can derive the following:

$$\begin{aligned}
\sum_{v \in V_i(t+1)} d_{G(t+1)}(v) &= \sum_{v \in X_0 \cap V_i(t+1)} d_{G(t+1)}(v) + \sum_{v \in V_i(t+1) \setminus X_0} d_{G(t+1)}(v) \le \mu |X_0'| + \left(1 - p + \frac{\varepsilon p}{5s}\right) \sum_{v \in V_i(t+1) \setminus X_0} d_{G(t)}(v) \\
&\le O(D^{19/10}) + \left(1 - p + \frac{\varepsilon p}{5s}\right)(1 + D^{-1/300}) |V_i(t+1)| D(t) \\
&\le \left(\left(1 - p + \frac{\varepsilon p}{5s}\right)(1 + D^{-1/300}) + D^{-1/300}\right) |V_i(t+1)| D(t) \\
&\le \left(1 - p + \frac{\varepsilon p}{4s}\right) |V_i(t+1)| D(t) - \left(\frac{\varepsilon p}{20s} - 2D^{-1/300}\right) |V_i(t+1)| D(t) \\
&\le \left(1 - p + \frac{\varepsilon p}{4s}\right) |V_i(t+1)| D(t).
\end{aligned}$$

From this, we conclude that

$$\bar{d}_{G(t+1)}(V_i(t+1)) \le \left(1 - p + \frac{\varepsilon p}{4s}\right) D(t) = D(t+1),$$

and therefore, for each $i \in I(t+1)$,

**F2** *the bad event $C_i$ that* **G2** *does not hold for i occurs with probability at most*

$$O(D) \exp(-\log D \log \log D).$$

Now, we prove the inequality ($\spadesuit$). Since $\mu = \lceil D^{9/10} \rceil$ and $S(t) = \Theta(D)$, $q := \frac{\mu p s}{S(t)} = \Theta(p D^{-1/10})$. For each $v \in X_j$, by the definition of $X_j$ we have

$$p_v = 1 - \frac{d_{G(t)}(v) p s}{S(t)} = 1 - (j \pm 1)\mu \frac{ps}{S(t)} = 1 - jq \pm q = 1 - jq \pm D^{-1/10}.$$

Since $j = O(D/\mu) \le O(D^{1/10})$, $jq \le O(p)$ and thus

$$|V_i(t)| = \sum_j |X_j| \ge \sum_j (1 - jq)|X_j| \ge \sum_j (1 - \frac{1}{2})|X_j| = \frac{1}{2} \sum_j |X_j| = \frac{1}{2} |V_i(t)|.$$

This implies $\sum_j (1 - jq)|X_j| = \Theta(|V_i(t)|) \overset{\mathbf{A0}}{=} \Theta(S(t)) = \Theta(D)$.

**Claim 4.15.** *The following holds:*

$$\frac{\sum_j(1-jq)j\mu|X_j|}{\sum_j(1-jq)|X_j|} \le D(t).$$

**Proof.** Let $M_1 = \sum_j(1-jq)j\mu|X_j|$, $M_2 = \sum_j(1-jq)|X_j|$, $N_1 = \sum_j |X_j|j\mu$, and $N_2 = \sum_j |X_j|$. Now, we have

$$Q_1 := M_1 \cdot N_2 = \sum_{j,k}(1-jq)j\mu|X_j||X_k| = \sum_{j,k}(1-kq)k\mu|X_j||X_k|;$$
$$Q_2 := N_1 \cdot M_2 = \sum_{j,k} |X_j|j\mu(1-kq)|X_k| = \sum_{j,k} |X_k|k\mu(1-jq)|X_j|.$$

It follows

$$2Q_1 = \mu\sum_{j,k}\Big((1-jq)j|X_j||X_k| + (1-kq)k|X_j||X_k|\Big) = \mu\sum_{j,k}\Big((j+k-(j^2+k^2)q)|X_j||X_k|\Big);$$
$$2Q_2 = \mu\sum_{j,k}\Big(|X_j|j(1-kq)|X_k| + |X_k|k(1-jq)|X_j|\Big) = \mu\sum_{j,k}\Big((j+k-2jkq)|X_j||X_k|\Big),$$

and therefore $Q_1 \le Q_2$. This implies

$$\frac{M_1}{M_2} \le \frac{N_1}{N_2} = \frac{\sum_j |X_j|j\mu}{\sum_j |X_j|} \le \frac{\sum_{v\in V_i(t)} d_{G(t)}(v)}{|V_i(t)|} \le D(t),$$

as desired. □

For each $j$, we rewrite the statement **(a)** by substituting $p_v$:

$$\left|X'_j\right| = (1-jq)|X_j| \pm \left(D^{-1/10}|X_j| + D^{2/5}\log^3 D\sqrt{|X_j|}\right)$$

Let $S := \{j \mid |X_j| \neq 0\}$. The definition of $X_j$ also implies $|S| = O(D^{1/10})$. Since $\sum_j |X_j| = |V_i(t)| = \Theta(D)$ and

$$\sum_j \sqrt{|X_j|} = \sum_{j\in S}\sqrt{|X_j|} \le \sqrt{|S|\sum_{j\in S}|X_j|} = O(\sqrt{D^{1/10}\cdot D}) = O(D^{11/20})$$

by Cauchy-Schwarz Inequality, we conclude

$$\sum_j\left(D^{-1/10}|X_j| + D^{2/5}\log^3 D\sqrt{|X_j|}\right) = O(D^{19/20}\log^2 D) \le D^{-1/25}\sum_j(1-jq)|X_j|.$$

Therefore,

$$|V_i(t+1)| = \sum_j\left|X'_j\right| \ge (1-D^{-1/25})\sum_j(1-jq)|X_j|,$$

and

$$
\begin{aligned}
&\sum_{v\in V_i(t+1)} d_{G(t)}(v) \le \sum_j |X'_j|(j+1)\mu = \sum_j |X'_j| j\mu + \mu|V_i(t+1)| \\
&\le \sum_j (1-jq)(j+1)\mu|X_j| + \sum_j (j+1)\mu\left(D^{-1/10}|X_j| + D^{2/5}\log^3 D\sqrt{|X_j|}\right) \\
&\le \sum_j (1-jq)j\mu|X_j| + \mu\sum_j (1-jq)|X_j| + \mu O(D^{1/10})D^{-1/25}\sum_j (1-jq)|X_j| \\
&\le \sum_j (1-jq)j\mu|X_j| + (\mu + \mu O(D^{1/10})D^{-1/25})\sum_j |X_j| \le \sum_j (1-jq)j\mu|X_j| + (\mu + \mu O(D^{3/50}))|V_i(t)| \\
&\overset{\mathbf{A0}}{=} \sum_j (1-jq)j\mu|X_j| + (\mu + \mu O(D^{3/50}))S(t) \le \sum_j (1-jq)j\mu|X_j| + D^{99/50}.
\end{aligned}
$$

It follows

$$
\begin{aligned}
&\frac{\sum_{v\in V_i(t+1)} d_{G(t)}(v)}{|V_i(t+1)|} \le \frac{\sum_j (1-jq)j\mu|X_j| + D^{99/50}}{(1-D^{-1/25})\sum_j (1-jq)|X_j|} \le \frac{1}{1-D^{-1/25}}\frac{\sum_j (1-jq)j\mu|X_j|}{\sum_j (1-jq)|X_j|} + O(D^{99/100}) \\
&\overset{\text{Claim 4.15}}{\le} (1+2D^{-1/25})D(t) + O(D^{99/100}) \le (1+2D^{-1/25}+D^{-1/200})D(t) \le (1+D^{-1/300})D(t).
\end{aligned}
$$

This proves the inequality (♠), as promised.

We now complete the proof of Theorem 4.11 by establishing the following claim, which implies that **G1** and **G2** simultaneously hold with a positive probability.

**Claim 4.16.**

$$
\mathbb{P}\left[\bigcap_{i\in I(t+1)} (\overline{\mathcal{B}_i} \cap \overline{C_i})\right] > 0.
$$

**Proof.** According to **F1** and **F2**, for each event $\mathcal{E}_i := \mathcal{B}_i \cup C_i$, we have

$$
\mathbb{P}[\mathcal{E}_i] \le D^2 \exp(-\log D \log\log D).
$$

Note that a fixed such event $\mathcal{E}_i$ corresponds to an index $i \in I(t+1)$, which in turn corresponds to a block $V_i(t)$ at the $t$-th round.

Define a graph $\Gamma$ with a vertex set consisting of $\{V_i(t) \mid i \in I(t)\}$ such that two vertices $V_i(t)$ and $V_j(t)$ are adjacent in $\Gamma$ if $e_{G(t)}(V_i(t), V_j(t)) > 0$. Given that $|V_i(t)| = S(t) = \Theta(D)$ and $\Delta(G(t)) = O(D)$, it follows that the maximum degree $\Delta(\Gamma)$ of the graph $\Gamma$ is $O(D^2)$. It is clear that the occurrence of the event $\mathcal{E}_i$ solely depends on the random selections made concerning blocks $V_j$ that lie within a maximum distance of two from $V_i$ according to our randomized algorithm. Consequently, $\Gamma^4$, which is derived from $\Gamma$ by connecting vertex pairs that are at most four steps apart, serves as a dependency graph for the set of events $\{\mathcal{E}_i\}_{i\in I(t)}$. Since $\Delta(\Gamma^4) \le \Delta(\Gamma)^4 = O(D^8) \le D^9 - 1$ and

$$
e \cdot D^2 \exp(-\log D \log\log D) \cdot D^9 < D^{12-\log\log D} < 1,
$$

by Lovász Local Lemma (Lemma 2.3), we have

$$\mathbb{P}\left[\bigcap_{i\in I(t+1)} \overline{\mathcal{E}_i}\right] > 0,$$

as desired. □

**Theorem 4.17.** *There is a multipartite graph $G$ and a $sD$-thick partition $\mathcal{P}$ such that the maximum block average degree is $D$, the local degree is 1, and there is no independent transversal blow-up with size $s$ with respect to $\mathcal{P}$.*

**Proof.** Let $G$ be a multipartite graph with a partition $\mathcal{P} = \{V_1, V_2, \ldots, V_{D+1}\}$ such that for each $i \in [D + 1]$, the set $V_i$ contains $sD$ vertices labeled as $v^i_1, v^i_2, \ldots, v^i_{sD}$. In graph $G$, two vertices are adjacent if and only if they share the same subscript (indicating they are labeled with the same number) but have different superscripts (indicating they belong to different sets $V_i$).

As a consequence, for each $j \in [sD]$, the set $\{v^1_j, v^2_j, \ldots, v^{D+1}_j\}$ forms a complete graph $K_{D+1}$ in $G$. Additionally, with respect to $\mathcal{P}$, the maximum block average degree of $G$ is $D$, and the local degree is 1. Note that the partition $\mathcal{P}$ is $sD$-thick.

Now, suppose $T = \{S_1, S_2, \ldots, S_{D+1}\}$ is an independent transversal blow-up of size $s$ with respect to $\mathcal{P}$, where $S_i \subseteq V_i$ for each $i$. Since the union of all $S_i$ has size $s(D + 1)$, which is greater than $sD$, there must exist some $j \in [sD]$ and distinct $i_1, i_2 \in [D + 1]$ such that $v^{i_1}_j \in S_{i_1}$ and $v^{i_2}_j \in S_{i_2}$. By the adjacency rule in $G$, $v^{i_1}_j$ and $v^{i_2}_j$ are adjacent. This contradicts the independence of $T$ (since no edge should exist between $S_{i_1}$ and $S_{i_2}$ in an independent transversal blow-up). Therefore, no independent transversal blow-up of size $s$ with respect to $\mathcal{P}$ exists in $G$. □

It is worth mentioning that by employing the same construction utilized in the aforementioned proof, we can derive the following theorem.

**Theorem 4.18.** *There is a multipartite graph $G$ and a $s\Delta$-thick partition $\mathcal{P}$ such that the maximum degree is $\Delta$, the local degree is 1, and there is no independent transversal blow-up with size $s$ with respect to $\mathcal{P}$.*

# 5 Many independent transversal blow-ups

**Lemma 5.1.** *Suppose that the following condition is met for some function $f(\Delta)$: if a graph $G$ has a maximum degree of at most $\Delta$ and possesses an $f(\Delta)$-thick partition $\mathcal{P}$, then each vertex in $G$ is included in some independent transversal with respect to $\mathcal{P}$.*

*Given this premise, if we consider another graph $G'$ that also has a maximum degree of at most $\Delta$ and features an $(f(\Delta)+\Delta)$-thick partition $\mathcal{P}'$, then $G'$ necessarily contains $f(\Delta)+\Delta$ vertex-disjoint independent transversals with respect to $\mathcal{P}'$.*

**Proof.** Let $\mathcal{P}' = (V_1', V_2', \ldots, V_r')$. Without loss of generality, assume that each block $V_i'$ has size exactly $\Delta_* := f(\Delta) + \Delta$. Now, it is sufficient to show that $G'$ has a factor of independent transversals respect to $\mathcal{P}'$.

Let $F := \{I_1, I_2, \ldots, I_{\Delta_*}\}$ such that

- $I_1, I_2, \ldots, I_{\Delta_*}$ are pair-wise disjoint, and each $I_i$ is an independent set of $G'$;
- $|I_i \cap V_j'| \leq 1$ for each $i \in [\Delta_*]$ and $j \in [r]$;
- $|\bigcup_{i=1}^{\Delta_*} I_i|$ (we call it the size of $F$) is as large as possible.

If every vertex of $G'$ is contained in some $I_i$, then we are done. Hence we assume, without loss of generality, that $V_1'$ has a vertex $v$ has not been covered by $F$. Since $|I_i \cap V_1'| \leq 1$ for each $i \in [\Delta_*]$ and $|V_1'| = \Delta_*$, we may further assume, without loss of generality, that $I_1 \cap V_1' = \emptyset$.

Let $J = \{j : |I_1 \cap V_j'| = 1\}$. For each $j \in J$, let $v_j$ be the vertex in $I_1 \cap V_j'$ and let $Z(j)$ collect all indices $i \in [\Delta_*]$ such that $v_j$ has a neighbor in $G'$ which is covered by $I_i$. For each $j \in [r]$, let $V_j$ be a subset of $V_j'$ obtained by removing all vertices which are contained in $I_i$ for every $i \in Z(j)$. Clearly, $|V_j| \geq |V_j'| - |Z(j)| \geq (f(\Delta) + \Delta) - \Delta = f(\Delta)$.

Let $G$ be the graph induced by $\bigcup_{i=1}^{r} V_i$. Note that $G$ has a maximum degree at most $\Delta$ and has a $f(\Delta)$-thick partition $\mathcal{P} := \{V_1, V_2, \ldots, V_r\}$. According to our assumption, the vertex $v$ is included in some independent transversal $T$ with respect to $\mathcal{P}$.

For each $j \in [r] \setminus \{1\}$, let $w_j \in T \cap V_j$. Note that $w_j$ may be $v_j$ itself. However, it is impossible that $w_j = v_j$ for all $j$ because otherwise we can add $v$ to $I_1$ without disturbing the first two conditions on $F$ but the size of $F$ is enlarged by 1.

For those $j$ such that $w_j \neq v_j$, if $w_j$ is contained in some $I_{c(j)}$, then we transfer $w_j$ to $I_1$ and move $v_j$ to $I_{c(j)}$; and if $w_j$ is not contained in any $I_i$, then we transfer $w_j$ to $I_1$ and relocate $v_j$ away from $I_1$. In each case, the updated set $I_{c(j)}$ is still independent because each neighbor of $v_j$ is not contained in $I_{c(j)}$ (otherwise $w_j$ had already been removed from $V_j'$ and we cannot see it in $V_j$). After this operation, we form a new $F$ with larger size because the new $I_1$ is exactly $T$, which owns at least one more vertex, say, $v$. This contradicts the maximality of $F$. □

**Lemma 5.2.** [1, Aharoni, Berger, and Ziv] *If $G$ is a multipartite graph with maximum degree at most $\Delta$ and a $2\Delta$-thick partition $\mathcal{P}$, then given any vertex $v$, there exists an independent transversal with respect to $\mathcal{P}$ containing $v$.*

**Lemma 5.3.** *If $G$ is a multipartite graph with maximum degree at most $\Delta$ and a $3\Delta$-thick partition $\mathcal{P}$, then $G$ has a factor of independent transversals with respect to $\mathcal{P}$.*

**Proof.** This is a consequence of Lemmas 5.1 and 5.2, resulting from substituting $f(\Delta)$ with $2\Delta$. □

**Theorem 5.4.** *Given a multipartite graph $G$ with a partition $\mathcal{P} = (V_1, V_2, \ldots, V_r)$ such that $|V_i| = 3s^2\Delta$ for each $i$, if $\Delta(G) \le \Delta$, then $G$ has a factor of independent transversal blow-ups with size $s$ with respect to $\mathcal{P}$.*

**Proof.** In this proof, we focus on the auxiliary graph $\ddot{G}$ as described in Definition 2, and our goal is to demonstrate that $\ddot{G}$ possesses a factor of independent transversals. According to Observation 2, $\ddot{G}$ exhibits a maximum degree of $\Delta_* := s\Delta$ and has a partition $\ddot{P}$ where each block has size exactly $3s\Delta = 3\Delta_*$. Consequently, by applying Lemma 5.3, we can conclude that $\ddot{G}$ indeed has a factor of independent transversals with respect to $\mathcal{P}$, fulfilling our objective. □

**Lemma 5.5.** [18, Loh and Sudakov] *For every $\varepsilon > 0$, the following holds. Given a multipartite graph $G$ with a partition $\mathcal{P} = (V_1, V_2, \ldots, V_r)$ such that $|V_i| = (2+\varepsilon)\Delta$ for each $i$, if $\Delta(G) \le \Delta$ and $\Gamma_{\mathcal{P}}(G) = o(\Delta)$, then $G$ has a factor of independent transversals with respect to $\mathcal{P}$.*

Here, we extend the result concerning independent transversals to include independent transversal blow-ups.

**Theorem 5.6.** *For every $\varepsilon > 0$, the following holds. Given a multipartite graph $G$ with a partition $\mathcal{P} = (V_1, V_2, \ldots, V_r)$ such that $|V_i| = (2+\varepsilon)s^2\Delta$ for each $i$, if $\Delta(G) \le \Delta$ and $\Gamma_{\mathcal{P}}(G) = o(\Delta)$, then $G$ has a factor of independent transversal blow-ups with size $s$ with respect to $\mathcal{P}$.*

**Proof.** Once again, we consider the auxiliary graph $\ddot{G}$ as outlined in Definition 2. According to Observation 2, $\ddot{G}$ possesses a maximum degree of $\Delta_* := s\Delta$ and a partition $\ddot{P}$ where each block has a precise size of $(2+\varepsilon)s\Delta = (2+\varepsilon)\Delta_*$. Additionally, with respect to $\ddot{P}$, $G$ has a local degree that is at most $so(\Delta) = o(\Delta_*)$. Consequently, Lemma 5.5 suggests that $\ddot{G}$ contains a factor of independent transversals with respect to $\ddot{P}$. Therefore, $G$ has a factor of independent transversal blow-ups, each of size $s$, with respect to $\mathcal{P}$. □

**Lemma 5.7.** [26, Wanless and Wood] *Fix integers $k \ge 2$ and $t \ge 1$. For a $k$-uniform hypergraph $G$, let $V_1, V_2, \ldots, V_r$ be a partition of $V(G)$ such that $|V_i| \ge t$ and at most $k^{-k}(k-1)^{k-1}t^{k-1}|V_i|$ stretched edges in $G$ intersect $V_i$ for each $i \in [n]$. Then, there exist at least $(\frac{k-1}{k}t)^r$ independent transversals of $V_1, V_2, \ldots, V_r$.*

**Theorem 5.8.** *Given a multipartite graph $G$ with a partition $\mathcal{P} = (V_1, V_2, \ldots, V_r)$, if for each $i$, $|V_i| \ge t$ and $V_i$ is incident with at most $\frac{t}{4s^2}|V_i|$ edges, then $G$ has at least*

$$\frac{1}{2^r}\binom{t}{s}^r$$

*independent transversal blow-ups with size $s$ with respect to $\mathcal{P}$. Consequently, if $\Gamma_{\mathcal{P}}(G) \le D$ and $\mathcal{P}$ is $4s^2D$-thick, then then $G$ possesses at least as many independent transversal blow-ups of size $s$ with respect to $\mathcal{P}$ as there are when substituting $4s^2D$ for $t$.*

**Proof.** For simplicity, we assume that the cardinality of each $V_i$ is $t$. We then consider the auxiliary graph $\dot{G}$ as described in Definition 1. According to Observation 1, $\dot{G}$ is a multipartite graph with a partition $\dot{\mathcal{P}} = (\dot{V}_1, \dot{V}_2, \ldots, \dot{V}_r)$ where, for each $i$, $|\dot{V}_i| = t_* := \binom{t}{s}$. Furthermore, $\dot{V}_i$ is incident with at most

$$s\binom{t-1}{s-1}\frac{t}{4s^2}|\dot{V}_i| = \frac{1}{4}t_*|\dot{V}_i|$$

edges. By applying Lemma 5.7, we deduce that $\dot{G}$ possesses $\left(\frac{t_*}{2}\right)^r$ independent transversals with respect to $\dot{\mathcal{P}}$. Consequently, $G$ has $\left(\frac{t_*}{2}\right)^r$ independent transversal blow-ups of size $s$ with respect to $\mathcal{P}$, as required. □

# 6 Remarks

So far, it may become apparent that the primary focus of our paper has been on proving Theorem 4.9. As previously discussed, if we substitute $s+\varepsilon$ in Theorem 4.9 with $s^2+\varepsilon$, then leveraging the auxiliary graph tools introduced in Section 3, the proof becomes significantly simpler, see Theorem 4.7. Furthermore, in this paper, we have derived a multitude of additional conclusions through the application of these auxiliary graph tools, all of which share a common condition featuring the expression $s^2$. To illustrate, refer to Theorems 4.3, 4.5, 5.4, 5.6, and 5.8 for concrete examples. Therefore, a natural question arises:

- do those theorems still hold when we substitute $s$ with $s^2$?

On the other hand, it would also be interesting to investigate the hypergraph version of our theorems.

# Acknowledgments

We express our gratitude to Professor Guanghui Wang for participating in our early discussions and offering many valuable suggestions. The third author also thanks Professor Guanghui Wang for the warm hospitality and assistance provided during his visit to Shandong University.